\documentclass[submission]{FPSAC2025}

\usepackage{xcolor}
\usepackage{multirow}
\usepackage{tikz}
\usepackage{ytableau}
\usepackage{mathrsfs}
\usepackage{thmtools}
\usepackage{thm-restate}
\usepackage{varwidth}
\usepackage{comment}
\usepackage{mathtools}
\usepackage[mathscr]{euscript}
\usepackage{bbm}
\usepackage{multicol}
\usepackage{enumerate}

\theoremstyle{definition}
\newtheorem{thm}{Theorem}[section]
\newtheorem{lem}[thm]{Lemma}
\newtheorem*{lem*}{Lemma}
\newtheorem*{thm*}{Theorem}
\newtheorem{prop}[thm]{Proposition}

\newtheorem{cor}[thm]{Corollary}

\newtheorem{defn}[thm]{Definition}
\newtheorem*{remark*}{Remark}
\newtheorem{remark}{Remark}
\newtheorem{example}{Example}

\newtheorem{cor/defn}[thm]{Corollary/Definition}
\DeclareMathOperator{\Par}{\mathbb{Y}}

\DeclareMathOperator{\sA}{\mathscr{A} }
\DeclareMathOperator{\sD}{\mathscr{D} }

\DeclareMathOperator{\Ind}{Ind}
\DeclareMathOperator{\Res}{Res}

\DeclareMathOperator{\Hilb}{\mathrm{Hilb}}

\DeclareMathOperator{\MOD}{\mathrm{Mod}}
\DeclareMathOperator{\bF}{\mathbb{F}}
\DeclareMathOperator{\Y}{\mathrm{Y}}
\DeclareMathOperator{\V}{\mathrm{V}}

\DeclareMathOperator{\ux}{\underline{x}}

\usepackage{lipsum}

\title[Higher Rank Macdonald Polynomials]{Higher Rank Macdonald Polynomials}

\author[Milo Bechtloff Weising]{Milo Bechtloff Weising\thanks{\href{mailto:milojbw@vt.edu}{milojbw@vt.edu}.}\addressmark{1}}

\address{\addressmark{1}Department of Mathematics, Virginia Tech, Blacksburg VA}

\received{\today}

\abstract{In this paper, we introduce higher rank generalizations of Macdonald polynomials. The higher rank non-symmetric Macdonald polynomials are Laurent polynomials in several sets of variables which form weight bases for higher rank polynomial representations of double affine Hecke algebras with respect to higher rank Cherednik operators. We prove that these polynomials satisfy generalized versions of the classical Knop--Sahi relations and we give combinatorial descriptions of their weights. The higher rank symmetric Macdonald polynomials are defined as Hecke-symmetrizations of the higher rank non-symmetric Macdonald polynomials and form eigenbases for the spaces of Hecke-invariant higher rank polynomials with respect to generalized finite variable Macdonald operators. We prove that the higher rank symmetric Macdonald polynomials satisfy stability properties allowing for the construction of infinite variable limits. These higher rank symmetric Macdonald functions form eigenbases for certain representations of the (positive) elliptic Hall algebra with respect to generalized infinite variable Macdonald operators. Lastly, we show that the higher rank polynomial representations may be used to construct higher rank polynomial representations of the double Dyck path algebra. This is a copy of the author's accepted extended abstract for FPSAC2025.}


\usepackage[backend=bibtex]{biblatex}
\begin{document}

\maketitle

\section{Introduction}

Macdonald polynomials are fundamental objects in algebraic combinatorics bridging combinatorics, representation theory, and geometry. In his highly influential paper \cite{MacDSLC}, Macdonald introduced the \textbf{\textit{symmetric Macdonald polynomials}} $P_{\lambda}(X;q,t)$ as $q,t$-analogues of the Schur functions $s_{\lambda}(x).$ In an effort to resolve open problems surrounding the Macdonald polynomials $P_{\lambda}$, Cherednik developed the theory of \textbf{\textit{double affine Hecke algebras}} (DAHAs) and \textbf{\textit{non-symmetric Macdonald polynomials}} \cite{C_95}. On the geometric side, Haiman showed \cite{haiman2000hilbert} that a plethystically transformed family of Macdonald polynomials $\tilde{H}(X;q,t)$ may be interpreted geometrically as the torus fixed point classes of the equivariant K-theory of the \textbf{\textit{Hilbert schemes}} $\Hilb_n(\mathbb{C}^2).$ Both of these algebraic and geometric perspectives have been developed and expanded over the last few decades leading to new areas in algebraic combinatorics.

Many authors have introduced and studied generalizations of Macdonald polynomials. Haiman introduced \textbf{\textit{wreath Macdonald polynomials}} which come from considering the K-theory of the fixed loci of $\Hilb_n(\mathbb{C}^2)$ with respect to certain cyclic group actions (see \cite{wreathmacdonald} for a survey). In another direction, Haglund introduced the \textbf{\textit{multi-t Macdonald polynomials}} which are combinatorially defined symmetric functions with infinitely many parameters $q,t_1,t_2,\ldots$ generalizing the usual modified Macdonald functions \cite{LOR_2023}. Dunkl--Luque defined \textbf{\textit{vector-valued}} generalizations of the non-symmetric and symmetric Macdonald polynomials \cite{DL_2011} which were extended by the author \cite{BWMurnaghanSLC} to give vector-valued generalizations of the infinite variable symmetric Macdonald functions. Gonz\'alez--Lapointe \cite{GL_2020} defined \textbf{\textit{multi-Macdonald polynomials}} which are related to the wreath Macdonald polynomials \cite{OS_private}. Other authors have studied \textbf{\textit{partially-symmetric}} generalizations of Macdonald polynomials and associated geometric constructions \cite{MBW} \cite{GCM_2017} \cite{goodberry2023} \cite{lapointe2022msymmetric} \cite{OBW_2024}.

In this paper, we take a different approach to constructing generalizations of the Macdonald polynomials. We start by introducing a representation-theoretic gadget called \textbf{\textit{expansion functors}} $\tau_{q_i,q_j}$ (Definition \ref{expansion functor def}) which build new DAHA modules from old ones and allow for varying choices of $q$-parameters in each copy of DAHA. Using these functors we construct the \textbf{\textit{higher rank polynomial representations}} of DAHA (Definition \ref{higher rank poly rep def}) by repeatedly applying these expansion functors to the ordinary polynomial representation of DAHA. These construction work in \textit{all} types but we focus only on type $\mathrm{GL}.$ The higher rank polynomial representations $V^{(r)}(q_1,\ldots,q_r)$ depend on the choice of parameters $q_1,\ldots,q_r$ and when they are generic the Cherednik elements $Y_1,\ldots,Y_n$ act with one-dimensional weight spaces. We define the \textit{\textbf{higher rank non-symmetric Macdonald polynomials}} $E_{\mu^{(1)},\ldots,\mu^{(r)}}(\ux_1,\ldots,\ux_{r};q_1,\ldots, q_r,t)$ as the $Y$-weight vectors of $V^{(r)}(q_1,\ldots,q_r)$ (with a specific choice for normalization) depending on $\mu^{(1)},\ldots,\mu^{(r)} \in \mathbb{Z}^n$ (Theorem \ref{main theorem 1}). In Proposition \ref{Higher Rank Knop-Sahi Relations} we derive higher rank analogues of the \textbf{\textit{Knop--Sahi relations}} for the usual non-symmetric Macdonald polynomials (Proposition \ref{usual Knop-Sahi}). We define the \textbf{\textit{higher rank symmetric Macdonald polynomials}} $P_{\nu^{\bullet}}$ (Definition \ref{sym macd def}) as certain Hecke-symmetrizations of the higher rank non-symmetric Macdonald polynomials and show they form a $Y_1+\ldots + Y_n$ eigenbasis for the \textbf{\textit{higher rank spherical representation}} $W_{n}^{(r)}(q_1,\ldots,q_r)$ with distinct weights (Theorem \ref{Sym Macd Thm}). In the last section, we consider stability for the higher rank symmetric Macdonald polynomials with respect to the natural quotient maps $x_{i,n+1} \rightarrow 0$ on the polynomial representations. We are able to define representations $W^{(r)}(q_1,\ldots,q_r)$ for the elliptic Hall algebra $\mathcal{E}^+$ along with \textit{\textbf{higher rank symmetric Macdonald functions}} $\mathcal{P}_{\nu^{\bullet}}$ which form an eigenbasis with respect to a higher rank Macdonald operator $\Delta$ (Theorem \ref{sym Macd function thm}). Furthermore, using machinery from the author's paper \cite{BWDDPA} we show that the spaces $W^{(r)}(q_1,\ldots,q_r)$ form the Hecke-invariants of a larger space $\mathcal{L}^{(r)}_{\bullet}(q_1,\ldots,q_r)$ 
called the \textit{\textbf{rank r polynomial representation}} of $\mathbb{B}_{t,q_1}$ (Corollary \ref{B tq cor}).

\section{Definitions and Notation}
\subsection{Combinatorics}

\begin{defn}
    Let $n \geq 1.$ We will write $\mathfrak{S}_n$ for the \textit{\textbf{symmetric group}} and $ \widehat{\mathfrak{S}}_n$ for the \textit{\textbf{extended affine symmetric group}}. The symmetric group $\mathfrak{S}_n$ has generators $s_1,\ldots, s_{n-1}$ given by simple adjacent transpositions $s_i:= (i,i+1)$ and $\widehat{\mathfrak{S}}_n$ has one additional generator $\pi$ satisfying $\pi^2 s_{n-1} = s_1\pi^2$ and $\pi s_i = s_{i+1}\pi$ for $1 \leq i \leq n-2.$ We will consider the natural action of $\widehat{\mathfrak{S}}_n$ on $\mathbb{Z}^n$ in which $\mathfrak{S}_n$ acts by permuting coordinates and $\pi(a_1,\ldots,a_n):= (a_n+1,a_1,\ldots,a_{n-1}).$ We will identify $\mathbb{Z}^n$ with the set $\widehat{\mathfrak{S}}_n/\mathfrak{S}_n$ of minimal length coset representatives and therefore consider $\mathbb{Z}^n \subset \widehat{\mathfrak{S}}_n.$
\end{defn}

\begin{defn}
Let $e_1,...,e_n$ be the standard basis of $\mathbb{Z}^n$ and let $\alpha \in \mathbb{Z}^n$. We define the\textbf{\textit{ Bruhat ordering}} on $\mathbb{Z}^n$, written simply by $<$, by first defining cover relations for the ordering and then taking their transitive closure. If $i<j$ such that $\alpha_i < \alpha_j$ then we say $\alpha > (ij)(\alpha)$ and additionally if $\alpha_j - \alpha_i > 1$ then $(ij)(\alpha) > \alpha + e_i - e_j$ where $(ij)$ denotes the transposition swapping $i$ and $j.$
\end{defn}

\subsection{Double Affine Hecke Algebras}

\begin{defn}
    Let $t,q_1, q_2,q_3,\ldots$ be a family of algebraically independent commuting free variables. Define $\mathbb{F}:= \mathbb{Q}(t,q_1,q_2,q_3,\ldots).$
\end{defn}

\begin{defn} \label{defn1}
Let $\ell \geq 1$.
Define the \textbf{\textit{double affine Hecke algebra}} $\mathscr{D}_n(q_{\ell})$ to be the $\bF$-algebra generated by $T_1,\ldots,T_{n-1}$, $X_1^{\pm 1},\ldots,X_{n}^{\pm 1}$, and $Y_1^{\pm 1},\ldots,Y_n^{\pm 1}$ with the following relations:

\begin{multicols}{2}
\begin{itemize}
    \item  
    \label{def-i} 
    $(T_i -1)(T_i +t) = 0$,
    \item [] $T_iT_{i+1}T_i = T_{i+1}T_iT_{i+1}$,
    \item [] $T_iT_j = T_jT_i$, $|i-j|>1$,
    \item 
    \label{def-ii} $T_i^{-1}X_iT_i^{-1} = t^{-1}X_{i+1}$,
    \item []$T_iX_j = X_jT_i$, $j \notin \{i,i+1\}$,
    \item []$X_iX_j = X_jX_i$,
    \item 
    \label{def-iii}$T_iY_iT_i = tY_{i+1}$,
    \item []$T_iY_j = Y_jT_i$, $j\notin \{i,i+1\}$,
    \item []$Y_iY_j = Y_jY_i$,
    \item 
    \label{def-iv}$Y_1T_1X_1 = X_2Y_1T_1$,
    \item 
    \label{def-v}$Y_1X_1\cdots X_n = q_{\ell} X_1\cdots X_nY_1$
    \item []
\end{itemize}
\end{multicols}

Define the special element $\pi$ by $$\pi := t^{-(n-1)}Y_1T_1\cdots T_{n-1}.$$ Define $\widetilde{\pi}:= X_1T_1^{-1}\cdots T_{n-1}^{-1}.$ The \textbf{\textit{affine Hecke algebra}} $\sA_n$ is the subalgebra of $\sD_n(q_{\ell})$ generated by the elements $T_1,\ldots,T_{n-1}, \pi, \pi^{-1}.$ We will write $\sD_n^{+}(q_{\ell})$ for the subalgebra of $\sD_n(q_{\ell})$ generated by $T_1,\ldots,T_{n-1},\pi,\widetilde{\pi}$ without containing the inverses of $\pi,\widetilde{\pi}.$ We write $\Y$ for the commutative subalgebra of $\sA_n$ generated by the \textit{Cherednik elements} $Y_1^{\pm 1},\ldots,Y_n^{\pm 1}.$ We may also consider the \textit{dual Cherednik elements} $\theta_1,\ldots,\theta_n$ defined by $\theta_i:= t^{i-1} T_{i-1}^{-1}\cdots T_{1}^{-1}\pi T_{n-1}\cdots T_{i}.$ We will write $\theta$ for the subalgebra of $\sA_n$ generated by $\theta_1^{\pm 1},\ldots,\theta_n^{\pm 1}$.

\end{defn}

\begin{remark}
    The algebra $\sA_n$ does not depend on the choice of parameter $q_{\ell}$ so we may consider $\sA_n$ as a subalgebra of \textit{every} $\sD_{n}(q_{\ell})$ for all $\ell \geq 1.$
\end{remark}

\begin{defn} \label{defn2}
Let $ \ell \geq 1$. The \textbf{\textit{ polynomial representation}} $\V(q_{\ell})$ of $\mathscr{D}_n(q_{\ell})$ is given by the following action on $\bF[x_1^{\pm 1},\ldots,x_n^{\pm 1}]$:

\begin{center}
\begin{itemize}
    \item $T_if(x_1,\ldots,x_n) = s_i f(x_1,\ldots,x_n) +(1-t)x_i \frac{1-s_i}{x_i-x_{i+1}}f(x_1,\ldots,x_n)$
    \item $X_if(x_1,..,x_n)= x_if(x_1,\ldots,x_n)$
    \item $\pi f(x_1,\ldots,x_n) = f(x_2,\ldots,x_{n},q_{\ell}^{-1}x_1).$
\end{itemize}
\end{center}

Here $s_i$ denotes the operator that swaps the variables $x_i$ and $x_{i+1}$.
\end{defn}

We will write $\epsilon^{(n)}$ for the normalized trivial Hecke idempotent which is the unique element of the \textbf{\textit{finite Hecke algebra}}, the subalgebra of $\sD_n(q_{\ell})$ generated by $T_1,\ldots,T_{n-1}$, satisfying $(\epsilon^{(n)})^2 = \epsilon^{(n)} \neq 0$ and for all $1 \leq i \leq n-1$, $T_i\epsilon^{(n)} = \epsilon^{(n)}T_i = \epsilon^{(n)}.$ We will write $\mathcal{S}_n^{+}(q_{\ell}):= \epsilon^{(n)} \sD_{n}^{+}(q_{\ell}) \epsilon^{(n)}$ for the \textit{\textbf{positive spherical double affine Hecke algebra}}.

\subsection{Non-Symmetric Macdonald Polynomials}
Fix some $\ell \geq 1.$ For the following definitions we will consider the polynomial representation of $\sD_n(q_{\ell}).$

\begin{defn}\label{non sym macd def rank 1}
The \textbf{\textit{non-symmetric Macdonald polynomials}} (for $GL_n$) are the family of Laurent polynomials $E_{\mu}(x_1,\ldots,x_n;q_{\ell},t) \in \bF[x_1^{\pm 1},\ldots,x_n^{\pm 1}]$ for $\mu \in \mathbb{Z}^n$ uniquely determined by the following:

\begin{itemize}
    \item $E_{\mu}$ has a triangular monomial expansion of the form $E_{\mu} = x^{\mu} + \sum_{\lambda < \mu} a_{\lambda}x^{\lambda}.$

    \item $E_{\mu}$ is a $\Y$-weight vector.
\end{itemize}
\end{defn}

For $\mu \in \mathbb{Z}^n$ the $\Y$-weight of $E_{\mu}$ is given by $\kappa_{\mu }= (\kappa_{\mu}(1),\ldots, \kappa_{\mu}(n))$ where 
$\kappa_{\mu}(j) := q_{\ell}^{-\mu_j}t^{\beta_{\mu}(j)}$ where 
$\beta_{\mu}(j) := \#(k < j |\mu_{k} > \mu_j) + \#(j< k |  \mu_j \leq \mu_k).$ Note that $\kappa_{\mu}(j)$ depends implicitly on the choice of $ \ell \geq 1.$

\begin{prop}[Knop--Sahi Relations]\label{usual Knop-Sahi}
    If $\mu \in \mathbb{Z}^n$ and $1\leq j \leq n-1$ such that $s_j(\mu) > \mu$, then 
    \begin{itemize}
        \item $E_{\pi(\mu)}(x_1,\ldots,x_n;q_{\ell},t) = q_{\ell}^{\mu_n}X_1\pi E_{\mu}(x_1,\ldots,x_n;q_{\ell},t)$
        \item $E_{s_j(\mu)}(x_1,\ldots,x_n;q_{\ell},t) = \left(T_j +\frac{t-1}{1-\frac{\kappa_{\mu}(j)}{\kappa_{\mu}(j+1)}}\right)E_{\mu}(x_1,\ldots,x_n;q_{\ell},t).$
    \end{itemize}
\end{prop}

\begin{remark}
    We may also define polynomials $F_{\mu}$ similar to the non-symmetric Macdonald polynomials which are $\theta$-weight vectors instead of $Y$-weight vectors. Much of the theory of non-symmetric Macdonald polynomials may be constructed for the $F_{\mu}$ polynomials with slightly different combinatorics. 
\end{remark}

\section{Representation Theory}
\subsection{Expansion Functors}

The constructions in this paper are centered around the following representation theoretic gadget.

\begin{defn}\label{expansion functor def}
    For any $i,j \geq 1$ define the \textit{\textbf{expansion functor}} $\tau_{q_i,q_j}: \sD_n(q_{j})-\MOD \rightarrow \sD_n(q_{i})-\MOD $ by 
    $$\tau_{q_i,q_j}(\V):= \Ind_{\sA_n}^{\sD_n(q_i)}\Res_{\sA_n}^{\sD_n(q_j)}(\V).$$
\end{defn}

These expansion functors make sense in \textit{all} Lie types and therefore so will most of the constructions in this paper. However, for the sake of clarity and brevity we will only focus on type $\mathrm{GL}_n.$ The following definition will be convenient for our purposes in this section.

\begin{defn}
    Let $\ell \geq 1.$ A $\sD_{n}(q_{\ell})$-module $\V$ is \textbf{\textit{exceptional}} if $\V$ is $\Y$-semisimple with \textit{distinct} $\Y$-weights  $\alpha = (\alpha(1),\ldots,\alpha(n))$ with $\alpha(i) \in t^{\mathbb{Z}}q_{\ell}^{\mathbb{Z}} q_{\ell+1}^{\mathbb{Z}}\cdots .$
\end{defn}

Crucially, the polynomial representations $\V(q_{\ell})$ are exceptional.

\begin{defn}
    For $\ell \geq 1$ and $\sigma \in \widehat{\mathfrak{S}}_n$ we define $\Psi_{q_{\ell}}^{\sigma}: (\mathbb{F}^{*})^{n} \rightarrow (\mathbb{F}^{*})^{n}$ by $\Psi_{q_{\ell}}^{1}(\alpha) = \alpha$, $\Psi_{q_{\ell}}^{\sigma \gamma} = \Psi_{q_{\ell}}^{\sigma}\Psi_{q_{\ell}}^{\gamma}$, and 
    \begin{multicols}{2}
    \begin{itemize}
        \item $\Psi_{q_{\ell}}^{s_i}(\alpha) = (\ldots, \alpha(i+1),\alpha(i), \ldots) $
        \item $\Psi_{q_{\ell}}^{\pi}(\alpha) = (q_{\ell}^{-1}\alpha(n),\alpha(1),\ldots).$
    \end{itemize}
    \end{multicols}
\end{defn}

Using Mackey Decomposition we find the following:

\begin{lem}\label{Mackey}
    If $\V$ is an exceptional $\sD_n(q_{i+1})$-module, then $\tau_{q_i,q_{i+1}}(\V)$ is an exceptional $\sD_n(q_{i})$-module. Further, if $\{\alpha_t | t \in I \}$ are the $\Y$-weights of $\V$, then $\{\Psi_{q_{i+1}}^{\mu}(\alpha_t) | (t,\mu) \in I \times \mathbb{Z}^n \}$ are the $\Y$-weights of $\tau_{q_{i},q_{i+1}}(\V)$.
\end{lem}

\subsection{Higher Rank Polynomial Representations}

\begin{defn}
    For $ r\geq 1$ we define the space of \textit{\textbf{rank r Laurent polynomials}} by 
    $$\mathbb{F}[x_{1,1}^{\pm 1},\ldots x_{1,n}^{\pm 1},\ldots,x_{r,1}^{\pm 1},\ldots x_{r,n}^{\pm 1} ] = \mathbb{F}[\underline{x}_1,\ldots,\underline{x}_r ]$$ where we use the shorthand $\underline{x}_i = (x_{i,1}^{\pm 1},\ldots, x_{i,n}^{\pm 1}).$ Further, for $\alpha \in \mathbb{Z}^n$ we will write 
    $\underline{x}_{i}^{\alpha}:= x_{i,1}^{\alpha_1}\cdots x_{i,n}^{\alpha_n}.$ We will consider the natural $r$-dimensional grading on $\mathbb{F}[\underline{x}_1,\ldots,\underline{x}_r ]$ given by 
    $\deg(\underline{x}_1^{\alpha^{(1)}}\cdots \underline{x}_r^{\alpha^{(r)}}):= (|\alpha^{(1)}|,\ldots, |\alpha^{(r)}|)$ where $|\alpha^{(i)}|:= \alpha^{(i)}_1+\cdots + \alpha^{(i)}_n \in \mathbb{Z}.$ For every $1\leq i \leq r$ and $1\leq j \leq n-1$ we write $\xi_j:\bF[\underline{x}_i] \rightarrow \bF[\underline{x}_i]$ for the operator 
    $\xi_j(\underline{x}_i^{\alpha}):= x_{i,j}\frac{\underline{x}_i^{\alpha}-\underline{x}_i^{s_j(\alpha)}}{x_{i,j}-x_{i,j+1}}.$
\end{defn}

\begin{defn}\label{higher rank poly rep def}
    For $r \geq 1$, we define the \textbf{\textit{rank r polynomial representation}} of $\sD_n(q_1)$ as the module 
    $$\V^{(r)}(q_1,\ldots, q_r):= \tau_{q_1,q_2}\cdots \tau_{q_{r-1},q_r} \V(q_r).$$ We naturally identify $\V^{(r)}(q_1,\ldots, q_r)$ as a $\mathbb{F}$-vector space with the space of rank $r$ Laurent polynomials $\mathbb{F}[\underline{x}_1,\ldots,\underline{x}_r ]$.
\end{defn}

\begin{remark}
    We stress that even though $\V^{(r)}(q_1,\ldots, q_r)$ is a module for $\sD_n(q_1)$, the operators defining this representation will involve all of the parameters $q_1,\ldots,q_r.$
\end{remark}

By repeatedly applying Lemma \ref{Mackey} we find:

\begin{cor}
    $\V^{(r)}(q_1,\ldots,q_r)$ is an exceptional $\sD_n(q_{1})$-module. 
\end{cor}

\begin{lem}
    The generators $T_j,X_i,\pi$ of $\sD_n(q_1)$ act on $\V^{(r)}(q_1,\ldots, q_r)$ according to the following formulas:
    \begin{itemize}
        \item $T_j$ acts as $ s_j^{\otimes r} + (1-t)\sum_{k=0}^{r-1}s_j^{\otimes k}\otimes \xi_j \otimes 1^{\otimes(r-1-k)}.$
        \item $X_i(\underline{x}_1^{\mu^{(1)}}\cdots \underline{x}_r^{\mu^{(r)}}) = x_{1,i} \cdot  \underline{x}_1^{\mu^{(1)}}\cdots \underline{x}_r^{\mu^{(r)}}$
        \item $\pi(\underline{x}_1^{\mu^{(1)}}\cdots \underline{x}_r^{\mu^{(r)}}) = q_1^{-\mu_n^{(1)}}\cdots q_r^{-\mu_n^{(r)}} x_{1,1}^{\mu^{(1)}_n}x_{1,2}^{\mu^{(1)}_1}\cdots x_{1,n}^{\mu^{(1)}_{n-1}}\cdots x_{r,1}^{\mu^{(r)}_n}x_{r,2}^{\mu^{(r)}_1}\cdots x_{r,n}^{\mu^{(r)}_{n-1}}.$
    \end{itemize}
\end{lem}

\section{Higher Rank Non-Symmetric Macdonald Polynomials}

\subsection{Main Construction}

In this section, we construct the higher rank non-symmetric Macdonald polynomials.

\begin{defn}\label{weights def}
    For $\mu^{\bullet} = (\mu^{(1)},\ldots, \mu^{(r)}) \in (\mathbb{Z}^n)^{r}$ define $\alpha_{\mu^{\bullet}} \in (\mathbb{F}^{*})^{n}$ by 
    $$\alpha_{\mu^{\bullet}}:= \Psi_{q_1}^{\mu^{(1)}}\cdots \Psi_{q_r}^{\mu^{(r)}}(t^{n-1},\ldots, t,1).$$
\end{defn}

\begin{example}
    For $\mu^{\bullet} = (0,1,0|2,0,0|0,0,1)$ we have $\alpha_{(0,1,0|2,0,0|0,0,1)} = (q_2^{-2}q_3^{-1},q_1^{-1}t,t^2).$
\end{example}

\begin{defn}
We define the homomorphism $\overline{\cdot} : \widehat{\mathfrak{S}}_n \rightarrow \mathfrak{S}_n$ via $\overline{s_i}:= s_i$ for $1\leq i \leq n-1$ and $\overline{\pi}:= s_{n-1}\cdots s_1.$ For $\mu \in \mathbb{Z}^n \equiv \widehat{\mathfrak{S}}_n/ \mathfrak{S}_n$ we write $\sigma_{\mu}:= \overline{\mu}\in \mathfrak{S}_n$ where we have identified $\mu$ with its corresponding minimal length coset representative. For $\mu^{\bullet}= (\mu^{(1)},\ldots,\mu^{(r)}) \in (\mathbb{Z}^n)^r$ we define $\gamma_{\mu^{\bullet}} \in (\mathbb{Z}^n)^r$ as 
$$\gamma_{\mu^{\bullet}}:= (\mu^{(1)},\sigma_{\mu^{(1)}}(\mu^{(2)}),\ldots, \sigma_{\mu^{(1)}}\cdots \sigma_{\mu^{(r-1)}}(\mu^{(r)})).$$ Further, we will write $\sigma_{\mu^{\bullet}}:= \sigma_{\mu^{(1)}}\cdots \sigma_{\mu^{(r)}}.$
\end{defn}

\begin{remark}
Note that $\gamma$ defines a bijection between $(\mathbb{Z}^n)^r$ and itself. For all $\mu^{\bullet} \in (\mathbb{Z}^n)^r$ and $1\leq i \leq n$ we have that 
$$\alpha_{\mu^{\bullet}}(i) = q_1^{-\gamma_{\mu^{\bullet}}^{(1)}(i)}\cdots q_r^{-\gamma_{\mu^{\bullet}}^{(r)}(i)}t^{n-\sigma_{\mu^{\bullet}}(i)}.$$
\end{remark}

\begin{defn}
    For $\mu \in \mathbb{Z}^n$ define $T_{\mu} \in \sA_n$ inductively via $T_{(0,\ldots,0)} := 1$ and
    \begin{itemize}
    \begin{multicols}{2}
        \item $T_{s_j(\mu)} := T_jT_{\mu}$ if $s_j(\mu) > \mu$
        \item $T_{\pi(\mu)}:= \pi T_{\mu}.$
        \end{multicols}
    \end{itemize}
\end{defn}

By induction using DAHA intertwiners we can show the following result:

\begin{thm}[Main Theorem]\label{main theorem 1}
    There exists a \textit{unique} family of homogeneous higher rank Laurent polynomials $$\{ E_{\mu^{(1)},\ldots, \mu^{(r)}}(\underline{x}_1,\ldots, \underline{x}_r;q_1,\ldots, q_r,t)  |~r \geq 1~, ~(\mu^{(1)},\ldots, \mu^{(r)}) \in (\mathbb{Z}^n)^r \}$$ satisfying the following properties:
    \begin{itemize}
        \item For $r = 1$ and $\mu \in \mathbb{Z}^n$, $E_{\mu}(\underline{x}_1;q_1,t) = E_{\mu}(x_{1,1},\ldots,x_{1,n};q_1,t)$ is the usual non-symmetric Macdonald polynomial for $\mu$ (Definition \ref{non sym macd def rank 1}). 
        \item Each $E_{\mu^{(1)},\ldots, \mu^{(r)}}(\underline{x}_1,\ldots, \underline{x}_r;q_1,\ldots, q_r,t) \in \V^{(r)}(q_1,\ldots,q_r)$ is a $\Y$-weight vector with weight $\alpha_{\mu^{(1)},\ldots, \mu^{(r)}}$.
        \item For every $\mu \in \mathbb{Z}^n$ and $(\beta^{(1)},\ldots, \beta^{(r)}) \in (\mathbb{Z}^n)^r$, we have for some higher rank Laurent polynomials $g_{\gamma}$ depending on $\mu$ and $(\beta^{(1)},\ldots, \beta^{(r)})$ the \textit{triangular} expansion 
        \begin{align*}
            &E_{\mu,\beta^{(1)},\ldots, \beta^{(r)}}(\underline{x}_1,\ldots,\underline{x}_{r+1};q_1,\ldots,q_{r+1},t) \\
            &= \underline{x}_1^{\mu}T_{\mu} E_{\beta^{(1)},\ldots, \beta^{(r)}}(\underline{x}_2,\ldots,\underline{x}_{r+1};q_2,\ldots,q_{r+1},t) + \sum_{\gamma < \mu } \underline{x}_1^{\gamma} g_{\gamma}(\underline{x}_2,\ldots, \underline{x}_{r+1}).\\
        \end{align*}
        
    \end{itemize}
\end{thm}

\begin{example}
\begin{align*}
    &E_{(0,1,0),(2,1,0)}(\ux_1,\ux_2;q_1,q_2,t) = x_{1,2}x_{2,1}^2 x_{2,2} + (t-1) x_{1,2}x_{2,1}^2x_{2,3} + \left(\frac{1-t}{1-q_2t^2} \right) x_{1,2}x_{2,1}x_{2,2}x_{2,3} \\
    &~~~~+ \left(\frac{1-t}{1-q_1q_2^{-2}t^{-2}} \right)x_{1,1}x_{2,2}^2x_{2,3} 
    + \left(\frac{1-t}{1-q_1q_2^{-2}t^{-2}} \right) \left(\frac{1-t}{1-q_2t^2} \right) x_{1,1}x_{2,1}x_{2,2}x_{2,3}.\\
\end{align*}
    
\end{example}

\begin{remark}
    The author conjectures that there exists a combinatorial formula for the higher rank non-symmetric Macdonald polynomials generalizing the Haglund--Haiman--Loehr formula \cite{haglund2007combinatorial}. The \textit{compression} methods of Guo--Ram \cite{GR_22} might be useful in finding such a formula.
\end{remark}

\begin{defn}
    For $r \geq 1$ and $\mu^{\bullet} \in (\mathbb{Z}^{n})^r$ we call  $E_{\mu^{\bullet}}(\underline{x}_1,\ldots,\underline{x}_r;q_1,\ldots,q_r,t)$ the \textit{\textbf{rank $r$ non-symmetric Macdonald polynomial}} corresponding to $\mu^{\bullet}.$ As a shorthand, we will often write $E_{\mu^{\bullet}}:= E_{\mu^{\bullet}}(\underline{x}_1,\ldots,\underline{x}_r;q_1,\ldots,q_r,t).$
\end{defn}

\begin{cor}
    For all $r \geq 1$ the rank r non-symmetric Macdonald polynomials $E_{\mu^{\bullet}}$ for $\mu^{\bullet} \in (\mathbb{Z}^n)^r$ are a $\Y$-weight basis for $\V^{(r)}(q_1,\ldots, q_r)$ with \textit{distinct} weights $\alpha_{\mu^{\bullet}}.$ 
\end{cor}

\begin{remark}
    We have the equalities:
    \begin{itemize}
        \item $E_{0,\mu^{\bullet}}(\underline{x}_1,\ldots,\underline{x}_{r+1};q_1,\ldots,q_{r+1},t) = E_{\mu^{\bullet}}(\underline{x}_2,\ldots,\underline{x}_{r+1};q_2,\ldots,q_{r+1},t)$
        \item $E_{\mu^{\bullet},0}(\underline{x}_1,\ldots,\underline{x}_{r+1};q_1,\ldots,q_{r+1},t) = E_{\mu^{\bullet}}(\underline{x}_1,\ldots,\underline{x}_{r};q_1,\ldots,q_{r},t).$
    \end{itemize}
    
\end{remark}

\begin{remark}
    We may analogously define higher rank non-symmetric Macdonald polynomials $F_{\mu^{\bullet}}$ using the dual Cherednik operators $\theta_i$ instead of the usual Cherednik operators $Y_i.$ These polynomials yield a $\theta$-weight basis for the spaces $\V^{(r)}(q_1,\ldots,q_r)$ with distinct weights. 
\end{remark}

\subsection{Higher Rank Knop--Sahi Relations}

Using Theorem \ref{main theorem 1} we may inductively prove the following:

\begin{prop}[Higher Rank Knop--Sahi Relations]\label{Higher Rank Knop-Sahi Relations}
Let $(\mu^{(1)},\ldots,\mu^{(r)}) \in (\mathbb{Z}^n)^r$, $1\leq j \leq n-1$, and $ 1 \leq \ell \leq r-1 $ such that $s_j(\mu^{(i)}) = \mu^{(i)}$ for $1 \leq i \leq \ell -1$ and $s_j(\mu^{(\ell)}) > \mu^{(\ell)}$, The following hold:
    \begin{itemize}
        \item $E_{\pi(\mu^{(1)}),\mu^{(2)},\ldots,\mu^{(r)}} = q_1^{\mu_n^{(1)}}X_1\pi E_{\mu^{\bullet}} $
        \item $E_{\mu^{(1)},\ldots,\mu^{(\ell -1)},s_j(\mu^{(\ell)}), \mu^{(\ell+1)},\ldots,\mu^{(r)}} = \left( T_j + \frac{t-1}{1-\frac{\alpha_{\mu^{\bullet}}(j)}{\alpha_{\mu^{\bullet}}(j+1)}} \right) E_{\mu^{\bullet}}.$
    \end{itemize} 
     
\end{prop}

Let $\omega:= (1,\ldots, 1) \in \mathbb{Z}^n.$ There is another non-trivial recurrence relation satisfied by the $E_{\mu^{\bullet}}.$

\begin{cor}\label{shifting cor}
    For all $1 \leq j \leq r$, $c \in \mathbb{Z}$, and $(\mu^{(1)},\ldots,\mu^{(r)}) \in (\mathbb{Z}^n)^r$, 
    $$E_{\mu^{(1)},\ldots,\mu^{(j)}+c \omega,\ldots, \mu^{(r)}} = \underline{x}_j^{c \omega} E_{\mu^{(1)},\ldots,\mu^{(r)}}.$$
\end{cor}

\section{Higher Rank Symmetric Macdonald Polynomials}

In this section, we will investigate the Hecke-invariants of the spaces $\V^{(r)}(q_1,\ldots,q_r)$. We use the word \textit{symmetric} to be synonymous with \textit{Hecke-invariant} instead of the more traditional meaning as \textit{permutation-invariant}.

\begin{remark}
    The rank $2$ polynomial
    $(T_1+t)(\underline{x}_{1}^{e_1} \underline{x}_{2}^{e_1}) = \ux_1^{e_1}\ux_2^{e_1}+\ux_{1}^{e_2}\ux_2^{e_2}+(1-t)\ux_1^{e_2}\ux_2^{e_1}$ is Hecke-invariant but clearly not symmetric under the obvious permutation action.
\end{remark}

\begin{defn}
    We will denote by $\Phi_n^{(r)}$ the set of all $(\nu^{(1)},\ldots,\nu^{(r)}) \in (\mathbb{Z}^n)^r$ such that $\nu^{(1)}$ is weakly decreasing and for all $1 \leq j \leq r-1$ whenever $\nu^{(j)}_{i}= \nu^{(j)}_{i+1}$ we have $\nu^{(j+1)}_{i} \geq \nu^{(j+1)}_{i+1}.$ We write $\Phi_{n,+}^{(r)}:= \Phi_n^{(r)}\cap (\mathbb{Z}_{\geq 0}^n)^r.$ 
\end{defn}

The set $\Phi_n^{(r)}$ indexes the set of orbits of $(\mathbb{Z}^n)^r$ under the diagonal permutation action of $\mathfrak{S}_n.$ Note that $\Phi_n^{(r)}$ naturally contains the set $\Par_n^{r} $ of all $r$-tuples of partitions $(\lambda^{(1)},\ldots,\lambda^{(r)})$ where each $\lambda^{(i)}$ has at most $n$ parts. 

The well founded-ness of the next definition is a nontrivial result involving the interplay between the dual Cherednik operators $\theta_i$ and the usual Cherednik operators $Y_i$.

\begin{defn}\label{sym macd def}
For $\nu^{\bullet} \in \Phi_n^{(r)}$ we define the \textit{\textbf{rank r symmetric Macdonald polynomial}} $P_{\nu^{\bullet}}$ to be the \textit{unique} multiple of $\epsilon^{(n)}(E_{\nu^{\bullet}})$ of the form
    $P_{\nu^{\bullet}} = F_{\mu^{\bullet}} + \sum_{\beta^{\bullet} \neq 
 \mu^{\bullet}} c_{\beta^{\bullet}} F_{\beta^{\bullet}}$ where $\mu^{\bullet}$ is the \textit{unique} index such that $F_{\mu^{\bullet}}$ has $\theta$-weight equal to the $Y$-weight $\alpha_{\nu^{\bullet}}$ of $E_{\nu^{\bullet}}$, $c_{\beta^{\bullet}}$ are some scalars, and $\beta^{\bullet}$ range over the indices of the $\theta$-weight vectors, except for $\mu^{\bullet}$, contained in the \textit{irreducible} $\sA_n$-submodule of $V^{(r)}(q_1,\ldots,q_r)$ generated by $F_{\mu^{\bullet}}$.
\end{defn}

Every $P_{\nu^{\bullet}}$ for $\nu^{\bullet} \in \Phi_n^{(r)}$ can be written in the form $P_{\nu^{\bullet}} = \underline{x_1}^{c_1\omega}\cdots \underline{x_r}^{c_r\omega} P_{\beta^{\bullet}} $ for $\beta^{\bullet} \in \Phi_{n,+}^{(r)}$ and some $c_1,\ldots, c_r \in \mathbb{Z}.$

\begin{defn}
    We define the \textbf{\textit{rank r spherical polynomial representation}} $W_{n}^{(r)}(q_1,\ldots,q_r)$ as the $\mathcal{S}_n^{+}(q_1)$-module 
    $W_{n}^{(r)}(q_1,\ldots,q_r):= \epsilon^{(n)}\left( V_{n,+}^{(r)}(q_1,\ldots,q_r) \right)$ where 
    $V_{n,+}^{(r)}(q_1,\ldots,q_r)$ is the subspace of polynomials $ \bF[x_{1,1},\ldots,x_{1,n},\ldots, x_{r,1},\ldots, x_{r,n}].$
\end{defn}

We will write $\Delta_n:= \epsilon^{(n)}\left((Y_1+\ldots+Y_{n}) - (1+t+\ldots + t^{n-1})\right)\epsilon^{(n)}.$ As a direct result of Theorem \ref{main theorem 1} we find:

\begin{thm}\label{Sym Macd Thm}
    The set of higher rank symmetric Macdonald polynomials $\{P_{\nu^{\bullet}}\}_{\nu^{\bullet}\in \Phi_{n,+}^{(r)}}$ is a $\Delta_n$-weight basis for $W_{n}^{(r)}(q_1,\ldots,q_r)$ with distinct spectrum given by 
    $$\Delta_n P_{\nu^{\bullet}} = \left(\sum_{i = 1}^{n} (q_1^{-\gamma_{\nu^{\bullet}}^{(1)}(i)}\cdots q_r^{-\gamma_{\nu^{\bullet}}^{(r)}(i)}-1)t^{n-\sigma_{\nu^{\bullet}}(i)} \right) P_{\nu^{\bullet}}.$$
\end{thm}

\section{Stability and $\mathbb{B}_{t,q}$}

Here we investigate stability phenomena for the higher rank Macdonald polynomials. In this section, we will write $\pi_n$ for the element $\pi \in \sD_n(q_{\ell})$ to differentiate between $\pi_m \in \sD_m(q_{\ell})$ for $n \neq m.$

\begin{defn}
    Define the map $\Pi^{(n)}: V_{n+1,+}^{(r)}(q_1,\ldots,q_r) \rightarrow V_{n,+}^{(r)}(q_1,\ldots,q_r)$ by   
$$\Pi^{(n)}(\ux_1^{\alpha^{(1)}}\cdots \ux_r^{\alpha^{(r)}}) := 
    \begin{cases}
        \ux_1^{(\alpha^{(1)}_1,\ldots,\alpha^{(1)}_{n})}\cdots \ux_r^{(\alpha^{(r)}_1,\ldots,\alpha^{(r)}_{n})} & \alpha^{(1)}_{n+1} = \ldots = \alpha^{(r)}_{n+1} = 0 \\
        0 & \text{otherwise} .\\
    \end{cases} $$
\end{defn}

The quotient maps $\Pi^{(n)}$ satisfy the following exceptional properties:

\begin{lem}\label{quotient relations lemma}
    The following relations hold:
    \begin{itemize}
    \begin{multicols}{2}
        \item $\Pi^{(n)}T_j = T_j\Pi^{(n)}$ for $1 \leq j \leq n-1$
        \item $\Pi^{(n)}\theta_i = \theta_i \Pi^{(n)}$ for $1\leq i \leq n$
        \item $\Pi^{(n)}(\theta_{n+1}-t^{n}) = 0$
        \item $\Pi^{(n)} \Delta_{n+1} = \Delta_n \Pi^{(n)}.$
    \end{multicols}
    \end{itemize}
\end{lem}

Using the properties of the maps $\Pi^{(n)}$, we find that the $P_{\nu^{\bullet}}$ are stable.

\begin{thm}\label{stability for sym}
    For all $\nu^{\bullet} \in \Phi_{n+1,+}^{(r)}$, 
     $$\Pi^{(n)}(P_{\nu^{\bullet}}) = 
    \begin{cases}
        P_{(\nu^{(1)}_1,\ldots, \nu^{(1)}_{n}),\ldots ,(\nu^{(r)}_1,\ldots, \nu^{(r)}_{n})} & \nu^{(1)}_{n+1}=\ldots = \nu^{(r)}_{n+1} = 0 \\
        0 & \text{otherwise} .\\
    \end{cases}$$
\end{thm}

\begin{defn}
    We may naturally include sets as $\Phi_{n,+}^{(r)} \rightarrow \Phi_{n+1,+}^{(r)}$ and define $\Phi^{(r)} := \lim_{\rightarrow} \Phi^{(r)}_{n,+}$ as the directed union. We will represent elements $\nu^{\bullet} \in \Phi^{(r)}$ as $$\nu^{\bullet} = (\nu^{(1)}_1,\ldots,\nu^{(1)}_{n_1}|\ldots | \nu^{(r)}_1,\ldots,\nu^{(r)}_{n_r})$$ for $n_j \geq 0$ with $\nu^{(j)}_{n_j} \neq 0$ and write $\ell(\nu^{\bullet}) = \max\{n_1,\ldots,n_r\}.$ Given $\nu^{\bullet} \in \Phi^{(r)}$ and $n \geq \ell(\nu^{\bullet})$ we will write $\iota_{n}(\nu^{\bullet}):= (\nu^{(1)}_{1},\ldots,\nu^{(1)}_{n_1},0,\ldots,0|\ldots | \nu^{(1)}_{1},\ldots,\nu^{(1)}_{n_r},0,\ldots,0) \in \Phi^{(r)}_{n,+}$. Let $\gamma_{\nu^{\bullet}}:= \gamma_{\iota_{\ell(\nu^{\bullet})}(\nu^{\bullet})}$ and $\sigma_{\nu^{\bullet}}:= \sigma_{\iota_{\ell(\nu^{\bullet})}(\nu^{\bullet})}.$ Define $W^{(r)}(q_1,\ldots,q_r)$ as the $\mathbb{Z}_{\geq 0}^r$-graded stable limit of the spaces $W_{n}^{(r)}(q_1,\ldots,q_r)$ with respect to the maps $\Pi^{(n)}.$ 
\end{defn}

Using \ref{stability for sym} we may define $\Delta = \lim_{n} \Delta_n$ as an operator on $W^{(r)}(q_1,\ldots,q_r).$ We denote the \textbf{\textit{positive elliptic Hall algebra}} of Burban--Schiffmann \cite{BS} \cite{SV} by $\mathcal{E}^{+} = \lim_{\leftarrow} \mathcal{S}_n ^{+}(q_1)$.

\begin{thm}\label{sym Macd function thm}
    For all $r\geq 1,$ $W^{(r)}(q_1,\ldots,q_r)$ is a graded $\mathcal{E}^{+}$-module with simple $\Delta$-spectrum. The $\Delta$-weight vectors are given for $\nu^{\bullet} \in \Phi^{(r)}$ by $\mathcal{P}_{\nu^{\bullet}}:= \lim_{n} P_{\iota_n(\nu^{\bullet})}$ with 
    $$\Delta \mathcal{P}_{\nu^{\bullet}} = \left(\sum_{i = 1}^{\ell(\nu^{\bullet})} (q_1^{-\gamma_{\nu^{\bullet}}^{(1)}(i)}\cdots q_r^{-\gamma_{\nu^{\bullet}}^{(r)}(i)}-1)t^{\ell(\nu^{\bullet})-\sigma_{\nu^{\bullet}}(i)} \right)\mathcal{P}_{\nu^{\bullet}}.$$
\end{thm}

\begin{remark}
    We call the $\mathcal{P}_{\nu^{\bullet}} \in W^{(r)}(q_1,\ldots,q_r)$ the \textit{\textbf{rank r symmetric Macdonald functions}}. For $r = 1$ these agree with the symmetric Macdonald functions $P_{\lambda}(X;q,t).$ If $\nu^{\bullet} \in \mathbb{Y}^{r}$ then the $\Delta$-weight of $\mathcal{P}_{\nu^{\bullet}}$ is simply $\sum_{i}^{\infty}(q_{1}^{-\lambda^{(1)}_i}\cdots q_{r}^{-\lambda^{(r)}_{i}}-1)t^{i-1}.$
\end{remark}

The algebra $\mathbb{B}_{q,t}$ was introduced by Carlsson--Gorsky--Mellit \cite{GCM_2017} and is a close relative of the double Dyck path algebra $\mathbb{A}_{q,t}$ introduced by Carlsson--Mellit \cite{CM_2015} in their proof of the Shuffle Theorem. The algebra $\mathbb{B}_{q,t}$ has been extended to $\mathbb{B}_{q,t}^{\text{ext}}$ by Gonz\'alez--Gorsky--Simental \cite{gonzález2023calibrated} to include operators $\Delta_{p_{\ell}}$ which are related to the Macdonald operators from classical Macdonald theory. In \cite{BWDDPA} the author introduced the concept of \textit\textbf{{compatible sequences}} of DAHA representations. Every compatible sequence of DAHA representations yields a corresponding representation of $\mathbb{B}_{q,t}$. Given the conventions of this paper, we will be looking at $\mathbb{B}_{t,q_1}$ specifically. Using machinery from \cite{BWDDPA} along with Lemma \ref{quotient relations lemma} we readily obtain the following.

\begin{cor}\label{B tq cor}
    For all $r \geq 1$ the families $\mathcal{C}^{(r)}(q_1,\ldots,q_r):= \left( V^{(r)}_{n,+}(q_1,\ldots,q_r), \Pi^{(n)} \right)_{n \geq 1}$ form \textit{compatible sequences}. As a consequence, there exist corresponding $\mathbb{B}_{t,q_1}$-modules $\mathcal{L}^{(r)}_{\bullet}(q_1,\ldots,q_r):= \mathcal{L}_{\bullet}(\mathcal{C}^{(r)}(q_1,\ldots,q_r))$ with $\mathcal{L}^{(r)}_0(q_1,\ldots,q_r) = W^{(r)}(q_1,\ldots,q_r).$  
\end{cor}

\begin{remark}
    For $r \geq 1$ we call the module $\mathcal{L}^{(r)}_{\bullet}(q_1,\ldots,q_r)$ the \textit{\textbf{rank r polynomial representations}} of $\mathbb{B}_{t,q_1}$. Since $\mathcal{L}_{\bullet}$ is functorial we in fact obtain graded $\mathbb{B}_{t,q_1}$-module maps $\mathcal{L}^{(r)}_{\bullet}(q_1,\ldots,q_r) \rightarrow \mathcal{L}^{(r+1)}_{\bullet}(q_1,\ldots,q_r,q_{r+1})$ for all $r \geq 1.$ Hence, by taking direct limits we may define the \textbf{\textit{universal polynomial representation}} of $\mathbb{B}_{t,q_1}$,
    $\mathcal{U}_{\bullet}(q_1,q_2,\ldots):= \lim_{\rightarrow} \mathcal{L}^{(r)}_{\bullet}(q_1,\ldots,q_r).$ The author conjectures that the modules $\mathcal{L}^{(r)}_{\bullet}(q_1,\ldots,q_r)$ yield extended actions by $\mathbb{B}_{t,q_1}^{\text{ext}}$ which in turn extend to a $\mathbb{B}_{t,q_1}^{\text{ext}}$ action on $\mathcal{U}_{\bullet}(q_1,q_2,\ldots).$ 
\end{remark}

\printbibliography

@article {CM_2015,
    AUTHOR = {Carlsson, Erik and Mellit, Anton},
     TITLE = {A proof of the shuffle conjecture},
   JOURNAL = {J. Amer. Math. Soc.},
  FJOURNAL = {Journal of the American Mathematical Society},
    VOLUME = {31},
      YEAR = {2018},
    NUMBER = {3},
     PAGES = {661--697},
      ISSN = {0894-0347},
   MRCLASS = {05E10 (05E05 33D52)},
  MRNUMBER = {3787405},
MRREVIEWER = {Tanja Stojadinovi\'{c}},
       DOI = {10.1090/jams/893},
       URL = {https://doi.org/10.1090/jams/893},
}

@misc{BWDDPA,
      title={Double Dyck Path Algebra Representations From DAHA}, 
      author={Bechtloff Weising, Milo},
      year={2024},
      eprint={2402.02843},
      archivePrefix={arXiv},
      primaryClass={math.RT},
      url={https://arxiv.org/abs/2402.02843}, 
}

@article {GCM_2017,
    AUTHOR = {Carlsson, Erik and Gorsky, Eugene and Mellit, Anton},
     TITLE = {The {$\Bbb{A}_{q,t}$} algebra and parabolic flag {H}ilbert
              schemes},
   JOURNAL = {Math. Ann.},
  FJOURNAL = {Mathematische Annalen},
    VOLUME = {376},
      YEAR = {2020},
    NUMBER = {3-4},
     PAGES = {1303--1336},
      ISSN = {0025-5831},
   MRCLASS = {16W50 (05E16 14C05 16E35 19G38 33C52 57K18)},
  MRNUMBER = {4081116},
MRREVIEWER = {Primo\v{z} Moravec},
       DOI = {10.1007/s00208-019-01898-1},
       URL = {https://doi.org/10.1007/s00208-019-01898-1},
}

@misc{gonzález2023calibrated,
      title={Calibrated representations of the double {D}yck path algebra}, 
      author={Nicolle González and Eugene Gorsky and José Simental},
      year={2023},
      eprint={2311.17653},
      archivePrefix={arXiv},
      primaryClass={math.RT}
}

@article {DL_2011,
    AUTHOR = {Dunkl, C. F. and Luque, J.-G.},
     TITLE = {Vector valued {M}acdonald polynomials},
   JOURNAL = {S\'{e}m. Lothar. Combin.},
  FJOURNAL = {S\'{e}minaire Lotharingien de Combinatoire},
    VOLUME = {66},
      YEAR = {2011/12},
     PAGES = {Art. B66b, 68},
      ISSN = {1286-4889},
   MRCLASS = {05E10 (05E05 05E15 33D80)},
  MRNUMBER = {2971011},
MRREVIEWER = {Himmet\ Can},
}

@article {haglund2007combinatorial,
    AUTHOR = {Haglund, J. and Haiman, M. and Loehr, N.},
     TITLE = {A combinatorial formula for nonsymmetric {M}acdonald
              polynomials},
   JOURNAL = {Amer. J. Math.},
  FJOURNAL = {American Journal of Mathematics},
    VOLUME = {130},
      YEAR = {2008},
    NUMBER = {2},
     PAGES = {359--383},
      ISSN = {0002-9327},
   MRCLASS = {05E05},
  MRNUMBER = {2405160},
MRREVIEWER = {Pavlo Pylyavskyy},
       DOI = {10.1353/ajm.2008.0015},
       URL = {https://doi.org/10.1353/ajm.2008.0015},
}

@misc{lapointe2022msymmetric,
      title={$m$-Symmetric functions, non-symmetric {M}acdonald polynomials and positivity conjectures}, 
      author={Luc Lapointe},
      year={2022},
      eprint={2206.05177},
      archivePrefix={arXiv},
      primaryClass={math.CO}
}

@article {haiman2000hilbert,
    AUTHOR = {Haiman, Mark},
     TITLE = {{H}ilbert schemes, polygraphs and the {M}acdonald positivity
              conjecture},
   JOURNAL = {J. Amer. Math. Soc.},
  FJOURNAL = {Journal of the American Mathematical Society},
    VOLUME = {14},
      YEAR = {2001},
    NUMBER = {4},
     PAGES = {941--1006},
      ISSN = {0894-0347},
   MRCLASS = {14C05 (05E05 20C30 33D45)},
  MRNUMBER = {1839919},
MRREVIEWER = {Claudio Procesi},
       DOI = {10.1090/S0894-0347-01-00373-3},
       URL = {https://doi.org/10.1090/S0894-0347-01-00373-3},
}

@article {SV,
    AUTHOR = {Schiffmann, Olivier and Vasserot, Eric},
     TITLE = {The elliptic {H}all algebra and the {$K$}-theory of the
              {H}ilbert scheme of {$\Bbb A^2$}},
   JOURNAL = {Duke Math. J.},
  FJOURNAL = {Duke Mathematical Journal},
    VOLUME = {162},
      YEAR = {2013},
    NUMBER = {2},
     PAGES = {279--366},
      ISSN = {0012-7094},
   MRCLASS = {19E08 (14C05 14C35 14F05 20C08)},
  MRNUMBER = {3018956},
MRREVIEWER = {Jens Hornbostel},
       DOI = {10.1215/00127094-1961849},
       URL = {https://doi.org/10.1215/00127094-1961849},
}

@article {BS,
    AUTHOR = {Burban, Igor and Schiffmann, Olivier},
     TITLE = {On the {H}all algebra of an elliptic curve, {I}},
   JOURNAL = {Duke Math. J.},
  FJOURNAL = {Duke Mathematical Journal},
    VOLUME = {161},
      YEAR = {2012},
    NUMBER = {7},
     PAGES = {1171--1231},
      ISSN = {0012-7094,1547-7398},
   MRCLASS = {16T05 (14F05)},
  MRNUMBER = {2922373},
MRREVIEWER = {Shilin\ Yang},
       DOI = {10.1215/00127094-1593263},
       URL = {https://doi.org/10.1215/00127094-1593263},
}

@article {MBW,
    AUTHOR = {Bechtloff Weising, Milo},
     TITLE = {Stable-Limit Non-symmetric {M}acdonald Functions in Type A},
   JOURNAL = {S\'em. Lothar. Combin.},
  FJOURNAL = {},
    VOLUME = {},
      YEAR = {2023},
    NUMBER = {89B.56},
     PAGES = {},
      ISSN = {},
   MRCLASS = {},
  MRNUMBER = {},
MRREVIEWER = {},
       DOI = {},
       URL = {},
}

@article {BWMurnaghanSLC,
    AUTHOR = {Bechtloff Weising, Milo},
     TITLE = {Murnaghan-Type Representations of the Elliptic Hall Algebra},
   JOURNAL = {S\'em. Lothar. Combin.},
  FJOURNAL = {},
    VOLUME = {},
      YEAR = {2024},
    NUMBER = {91B.47},
     PAGES = {},
      ISSN = {},
   MRCLASS = {},
  MRNUMBER = {},
MRREVIEWER = {},
       DOI = {},
       URL = {},
}

@article {MacDSLC,
    AUTHOR = {Macdonald, Ian},
     TITLE = {A New Class of Symmetric Functions},
   JOURNAL = {S\'em. Lothar. Combin.},
  FJOURNAL = {},
    VOLUME = {},
      YEAR = {1988},
    NUMBER = {B20a},
     PAGES = {},
      ISSN = {},
   MRCLASS = {},
  MRNUMBER = {},
MRREVIEWER = {},
       DOI = {},
       URL = {},
}

@Article{C_95,
 Author = {Cherednik, Ivan},
 Title = {Double affine {Hecke} algebras and {Macdonald}'s conjectures},
 FJournal = {Annals of Mathematics. Second Series},
 Journal = {Ann. Math. (2)},
 ISSN = {0003-486X},
 Volume = {141},
 Number = {1},
 Pages = {191--216},
 Year = {1995},
 Language = {English},
 DOI = {10.2307/2118632},
 zbMATH = {750768},
 Zbl = {0822.33008}
}

@misc{OBW_2024,
    title={Stable-limit partially symmetric Macdonald functions and parabolic flag Hilbert schemes},
    author={Daniel Orr and Milo Bechtloff Weising},
    year={2024},
    eprint={2410.13642},
    archivePrefix={arXiv},
    primaryClass={math.CO}
}

@misc{goodberry2023,
      AUTHOR = {Goodberry, Ben},
      EPRINT = {2311.12216},
  EPRINTTYPE = {arxiv},
       TITLE = {Type A Partially-Symmetric Macdonald Polynomials},
        YEAR = {2023}}

@article{OS_private,
    AUTHOR = {Orr, Daniel and Shimozono, Mark},
     TITLE = {Private Communication},
   JOURNAL = {},
  FJOURNAL = {},
    VOLUME = {},
      YEAR = {2024},
    NUMBER = {},
     PAGES = {},
      ISSN = {},
   MRCLASS = {},
  MRNUMBER = {},
MRREVIEWER = {},
       DOI = {},
        eprint={},
       URL = {}
}

@Article{GR_22,
 Author = {Guo, Weiying and Ram, Arun},
 Title = {Comparing formulas for type {{\(GL_n\)}} {Macdonald} polynomials},
 FJournal = {Algebraic Combinatorics},
 Journal = {Algebr. Comb.},
 ISSN = {2589-5486},
 Volume = {5},
 Number = {5},
 Pages = {849--883},
 Year = {2022},
 Language = {English},
 DOI = {10.5802/alco.227},
 zbMATH = {7615644},
 Zbl = {1503.05122}
}

@Article{wreathmacdonald,
      title={Wreath Macdonald polynomials, a survey}, 
      author={Daniel Orr and Mark Shimozono},
      year={2023},
      eprint={2308.12166},
      archivePrefix={arXiv},
      primaryClass={math.CO},
      url={https://arxiv.org/abs/2308.12166}, 
}

@Article{LOR_2023,
 Author = {Lee, Seung Jin and Oh, Jaeseong and Rhoades, Brendon},
 Title = {Haglund's conjecture for multi-{{\(t\)}} {Macdonald} polynomials},
 FJournal = {Discrete Mathematics},
 Journal = {Discrete Math.},
 ISSN = {0012-365X},
 Volume = {346},
 Number = {6},
 Pages = {12},
 Note = {Id/No 113360},
 Year = {2023},
 Language = {English},
 DOI = {10.1016/j.disc.2023.113360},
 zbMATH = {7676440},
 Zbl = {1519.05248}
}

@Article{GL_2020,
 Author = {Gonz{\'a}lez, Camilo and Lapointe, Luc},
 Title = {Multi-{Macdonald} polynomials},
 FJournal = {Discrete Mathematics},
 Journal = {Discrete Math.},
 ISSN = {0012-365X},
 Volume = {343},
 Number = {12},
 Pages = {13},
 Note = {Id/No 112111},
 Year = {2020},
 Language = {English},
 DOI = {10.1016/j.disc.2020.112111},
 zbMATH = {7257933},
 Zbl = {1448.05208}
}

\end{document}